\documentclass[a4paper]{article}
\usepackage[all]{xy}
\usepackage[utf8]{inputenc}
\usepackage{amsmath,bm}
\usepackage{amssymb}
\usepackage{amsthm}
\usepackage{mathtools}
\usepackage{cases}
\usepackage{float}
\usepackage{here}
\usepackage{comment}
\usepackage{ulem}
\usepackage[T1]{fontenc}
\usepackage{textcomp}
\usepackage{graphicx, color}
\usepackage{tikz,epic,eso-pic}

\theoremstyle{plain}
\newtheorem{thm}{Theorem}[section]
\newtheorem{lem}[thm]{Lemma}
\newtheorem{cor}[thm]{Corollary}
\newtheorem{prop}[thm]{Proposition}

\theoremstyle{definition}
\newtheorem{df}[thm]{Definition}
\newtheorem{eg}[thm]{Example}

\theoremstyle{remark}
\newtheorem{rem}[thm]{Remark}

\numberwithin{equation}{section}
\allowdisplaybreaks[3]

\newcommand{\Z}{\mathbb{Z}}

\newcommand{\del}{\partial}
\newcommand{\wh}{\widehat}
\newcommand{\too}{\longrightarrow}

\DeclareMathOperator{\MC}{MC}
\DeclareMathOperator{\MH}{MH}

\DeclareMathOperator{\rk}{rk}

\DeclareMathOperator{\dMH}{\mathcal{MH}}

\title{
Magnitude homology and Path homology
}

\author{Yasuhiko \textsc{Asao}\thanks{Department of Applied Mathematics, Fukuoka University \texttt{asao@fukuoka-u.ac.jp}}}
\date{\today}

\begin{document}
\maketitle
\begin{abstract}
In this article, we show that \textit{magnitude homology} and \textit{path homology} are closely related, and we give some applications. We define differentials $\MH^{\ell}_k(G) \too \MH^{\ell-1}_{k-1}(G)$ between magnitude homologies of a digraph $G$, which make them chain complexes. Then we show that its homology $\dMH^{\ell}_k(G)$ is non-trivial and  homotopy invariant in the context of `homotopy theory of digraphs' developed by Grigor'yan--Muranov--S.-T. Yau \textit{et al} (G-M-Ys in the following). It is remarkable that the diagonal part of our homology $\dMH^{k}_k(G)$ is isomorphic to the \textit{reduced path homology} $\tilde{H}_k(G)$ also introduced by G-M-Ys. Further, we construct a spectral sequence whose first page is isomorphic to magnitude homology  $\MH^{\ell}_k(G)$, and the second page is isomorphic to our homology $\dMH^{\ell}_k(G)$. As an application, we show that the diagonality of magnitude homology implies triviality of reduced path homology. We also show that $\tilde{H}_k(g) = 0$ for $k \geq 2$ and $\tilde{H}_1(g) \neq 0$ if any edges of an undirected graph $g$ is contained in a cycle of length $\geq 5$.
\end{abstract}
\section{Introduction}
\textit{Magnitude homology} of a graph $ \MH^{\ell}_k(g)$ is first introduced by Hepworth-Willerton in \cite{HW} as a categorification of magnitude, a categorical analogue of Euler characteristic introduced by Leinster (\cite{L1}). It is shown by Gu (\cite{Gu} Appendix A) that such a categorification truly contains more information than magnitude. Magnitude homology can be defined for general metric spaces (\cite{LS}), and several authors show interesting geometric properties of that. It is remarkable that such properties are similar to those of singular homology theory (\cite{A}, \cite{AHK}, \cite{Gom}, \cite{HW}, \cite{YK}, \cite{TY}, for example). However, only a few relations to the other invariants have been shown so far. For example, the \textit{girth} of a graph controls magnitude homology to some extent (\cite{AHK}). For geodesic metric spaces, magnitude homology indicates curvature boundedness or uniqueness of geodesics (\cite{A}, \cite{Gom}).

On the other hand, A. Grigor'yan, R. Jimenez, Y. Lin, Y. Muranov, and S.-T. Yau have introduced and developed the theory of \textit{path homology} of \textit{directed graphs (digraphs)} in a large literature of their works (\cite{GJMY}, \cite{GJMY2}, \cite{GLMY}, \cite{GMY}, for example). While there has been a number of attempts to define a (co)homology theory for (di)graphs, their motivation is to find a good one which can be non-trivial in all dimensions, and satisfies suitable axioms of homology theory. Indeed, path homology possesses homotopy invariance in a sense, and can be considered as a pivotal object to develop a homotopy theory for digraphs (\cite{GLMY}). Further, it is shown that path homology satisfies an analogue of the Eilenberg--Steenrod axioms (\cite{GJMY}). One of the remarkable applications of their theory is that they give a new elementary proof of a theorem of Gerstenhaber--Schack that identifies simplicial cohomology as a Hochschild cohomology (\cite{GMY}). There are also several studies considering path homology of hypergraphs, multigraphs and quivers, involving authors other than the above researchers (\cite{GJMY2}, \cite{GMVY}). We denote the $k$-th reduced path homology of a digraph $G$ by $\tilde{H}_k(G)$.

We remark here that the category of graphs and graph maps can be embedded to the category of digraphs and digraph maps (\textit{Remark} \ref{embed}). We denote this embedding functor by $g \mapsto \overrightarrow{g}$ for an undirected graph $g$. It is known that the definition of magnitude homology of undirected graphs can be extended to digraphs since they can also be metrized by directed path metric. We employ this extension throughout this article. We denote a digraph with finite vertices by $G$, and an undirected graph with finite vertices by $g$.

In this article, we show that magnitude homology and path homology are closely related, and we give some applications. We construct homomorphisms $\del'^{\ell}_k : \MH^{\ell}_k(G) \too \MH^{\ell-1}_{k-1}(G)$ satisfying $\del'^{\ell-1}_{k-1}\circ \del'^{\ell}_k$, and we denote its homology by $\dMH^{\ell}_k(G)$ (Definition \ref{defdmh}). Then we have the following.
\begin{thm}[Proposition \ref{htpyequiv}]
If two digraphs $G, G'$ are homotopy equivalent, then we have $\dMH^{\ell}_k(G) \cong \dMH^{\ell}_k(G')$ for every $\ell, k \in \Z$.
\end{thm}
\begin{thm}[Proposition \ref{pathdmh}]
We have $\tilde{H}_k(G) \cong \dMH^{k}_{k}(G)$ for every digraph $G$ and every $k \in \Z$.
\end{thm}
Here we employ the notion of homotopy for digraphs (Definition \ref{defhtpy}) introduced in \cite{GJMY} and \cite{GLMY}. We also show non-triviality of $\dMH^{\ell}_k(G)$ for $\ell \neq k$ in Proposition \ref{nontrivial}.

Further, we construct a spectral sequence whose first page is isomorphic to magnitude homology $\MH^{\ell}_k(G)$, and the second page is isomorphic to our homology $\dMH^{\ell}_k(G)$ (Proposition \ref{magspec}). The same spectral sequence is just mentioned in Remark 8.7 of \cite{HW}. As applications, we have the following. We say that a digraph $G$ is \textit{diagonal} if $\MH^{\ell}_k(G) = 0$ for $\ell \neq k$. 
\begin{thm}[Proposition \ref{diagvanish}]\label{thm3}
 If a digraph $G$ is diagonal and has finite diameter, then the reduced path homology $\tilde{H}_k(G)$ is $0$ for every $k \in \Z$.
\end{thm}
Note that $\overrightarrow{g}$ has finite diameter whenever an undirected graph $g$ is connected. We show that the inverse of Theorem \ref{thm3} is not true in Example \ref{counter}. We apply a result in \cite{AHK} to show the following. See also \textit{Remark} \ref{errortum}. We define the \textit{girth} of an undirected graph $g$ by the minimum length of cycles in $g$.
\begin{thm}[Proposition \ref{girth5}]
 If the girth of an undirected graph $g$ is $\geq 5$, then we have $\tilde{H}_k(\overrightarrow{g}) = 0$ for $k\geq 2$.
\end{thm}
\begin{thm}[Proposition \ref{girth52}]
  If the girth of an undirected graph $g$ is $ \geq 5$ and $< \infty$, then  we have $\tilde{H}_1(\overrightarrow{g}) \neq 0$.
\end{thm}
The rest of this article is organized as follows. In section \ref{sectdef}, we define chain complexes we consider throughout this article, in particular we define differentials on magnitude homologies. In section \ref{funct}, we show functoriality of the homologies defined in the previous section. In section \ref{sectop}, we consider some operations on the chain complexes defined in the previous section. These operations are used to prove the homotopy invariance in the latter section. In section \ref{sectinv}, we show the homotopy invariance of our homologies. In section \ref{sectcoin}, we show that our homology $\dMH^{k}_k(G)$ coincides with the reduced path homology $\tilde{H}_k(G)$. In section \ref{sectspec}, we construct a spectral sequence mentioned above after a short introduction to this subject. In section \ref{sectapl},  we give some applications to magnitude homology and path homology.
\subsubsection*{Acknowledgement}
The author thanks to M. Yoshinaga for helpful private communication, especially for pointing out \textit{Remark} \ref{yoshinaga}. He also expresses gratitude to an anonymous referee for reading the script in detail and giving him very helpful suggestions.
\section{Definitions of chain complexes}\label{sectdef}
In this section, we define chain complexes $\Lambda_\ast(G), R_\ast(G), \MC^{\ell}_\ast(G)$ and $\MH^{\ast}_\ast(G)$ that we consider throughout this article. Note that a part of our definition overlaps with that of magnitude homology introduced in \cite{HW} and \cite{LS}. Although they can be defined in a more concise way by using the spectral sequence introduced in Section \ref{sectspec}, we define them `by hand' to make it as accessible as possible to readers not familiar with that subject.
\begin{df}
A \textit{directed graph or digraph} $G$ is a pair of sets $(V(G), E(G))$ such that $E(G) \subset \{(x, y) \in V(G) \times V(G) \mid x \neq y\}$. We call elements of $V(G)$ and $E(G)$ \textit{vertices} and \textit{edges} respectively. A \textit{directed path} between  vertices $x$ and $y$ in $G$ is a sequence of edges $(x_0, y_0), \dots, (x_n, y_n)$ with $x_0 = x,  y_n = y,$ and $y_i = x_{i+1}$ for $0 \leq i \leq n-1$. We define the length of such a sequence by $n$. We also define the  \textit{directed path metric} $d$ on $G$ by $d(x, y) = $ `minimum length of directed paths between $x$ and $y$' if there exists a directed path between $x$ and $y$, and $d(x, y) = \infty$ if there does not.
\end{df}
Throughout this article, let $G$ be a digraph with $\#V(G) < \infty$. Let $d$ be the directed path metric on $G$.
\begin{df}
We define modules $\Lambda_k(G)$ for $k \in \Z$ as follows. For $k \geq -1$, we define  $
\Lambda_k(G) = \Z\langle (x_0, \dots, x_k) \in V(G)^{k+1} \rangle$, which is a free $\Z$ module generated by $(k+1)$-tuples of vertices of $G$, admitting the unique $0$-tuple $()$. We define $\Lambda_{k}(G) = 0$ for $k \leq -2$. 
\end{df}
\begin{df}
 For a $(k+1)$-tuple $ (x_0, \dots, x_k) \in V(G)^{k+1}$, we define 
 \[
 L(x_0, \dots, x_k) = \begin{cases} \sum_{i = 0}^{k-1} d(x_i, x_{i+1})  & k \geq 0, \\ -1 & k = -1.\end{cases}
 \]
 
\end{df}
\begin{df}
For $k, \ell \in \Z$, we define a homomorphism 
\[
(-)^{\leq \ell} : \Lambda_k(G) \too \Lambda_k(G)
\]
by 
\[
(x_0, \dots, x_k)^{\leq \ell} = \begin{cases} (x_0, \dots, x_k) & \text{if $L(x_0, \dots, x_k) \leq \ell$} \\ 0 & \text{otherwise}.\end{cases}
\]
We denote homomorphism $\mathrm{Id}_{\Lambda_{k}(G)} - (-)^{\leq \ell}$ and $(-)^{\leq \ell} - (-)^{\leq \ell - 1}$ by $(-)^{\geq \ell + 1}$ and $(-)^{\ell}$ respectively.
\end{df}
\begin{df}
For $k \in \Z$, we define a homomorphism $\tilde{\del}_k : \Lambda_k(G) \too \Lambda_{k-1}(G)$ by 
\[
\tilde{\del}_k (x_0, \dots, x_k) = \sum_{i = 0}^{k} (-1)^i (x_0, \dots, \check{x}_i, \dots, x_k).
\]
Here $\check{x}_i$  means omitting the $i$-th term.
\end{df}
The following can be easily checked.
\begin{lem}\label{id0}
For $k, \ell, \ell' \in \Z$ with $\ell \leq \ell'$, we have $(\tilde{\del}_k x^{\leq \ell})^{\leq \ell'} = \tilde{\del}_k x^{\leq \ell}$ and $(\tilde{\del}_k x^{\leq \ell})^{\geq \ell' + 1} = 0$ for every $x \in \Lambda_k(G)$.
\end{lem}
\begin{lem}\label{lambda}
For every $k \in \Z$, $\tilde{\del}_{k-1}\circ \tilde{\del}_k = 0$. Hence $(\Lambda_\ast, \tilde{\del}_\ast)$ is a chain complex.
\end{lem}
\begin{df}
We define submodules $I_k(G)$ of $\Lambda_k(G)$ for $k \in \Z$ as follows. For $k \geq 1$, we define  
\[
I_k(G) = \Z\langle (x_0, \dots, x_k) \in V(G)^{k+1} \mid x_i = x_{i+1}\ \text{for some $0 \leq i \leq k-1$}\rangle.\]
We define $I_{k}(G) = 0$ for $k \leq 0$. 
\end{df}
It is easily checked that we can restrict $\tilde{\del}_k$'s to $I_k(G) \too I_{k-1}(G)$.
\begin{df}
We define modules $R_k(G)$ for $k \in \Z$ by $R_k(G) = \Lambda_k(G)/I_k(G)$. Then the quotient homomorphism of $\tilde{\del}_k$  is well-defined, and we denote it by $\del_k : R_k(G) \too R_{k-1}(G)$.
\end{df}
The following are clear from Lemmas \ref{id0} and  \ref{lambda}.
\begin{lem}\label{id}
For $k, \ell, \ell' \in \Z$ with $\ell \leq \ell'$, we have $(\del_k x^{\leq \ell})^{\leq \ell'} = \del_k x^{\leq \ell}$ and $(\del_k x^{\leq \ell})^{\geq \ell' + 1} = 0$ for every $x \in R_k(G)$.
\end{lem}
\begin{lem}\label{rchain}
$(R_\ast, \del_\ast)$ is a chain complex.
\end{lem}
For $k \geq 1$, the module $R_k(G)$ is isomorphic to a submodule of $\Lambda_k(G)$ defined by
\[
 \Z\langle (x_0, \dots, x_k) \in V(G)^{k+1} \mid x_i \neq x_{i+1}\ \text{for any $0 \leq i \leq k-1$}\rangle.\]
 We also have $R_k(G) \cong \Lambda_k(G)$ for $k \leq 0$. Hence homomorphisms $(-)^{\leq \ell}, (-)^{\geq \ell}$, and $(-)^{\ell}$ are well-defined on $R_k(G)$'s. We denote it by the same symbols. In the following, we consider modules $R_k(G)$'s and homomorphisms defined on them.
\begin{lem}\label{easylemma}
For $\ell, k \in \Z$, we have
\[
\del_{k-1}(\del_k x)^{\ell} = -\del_{k-1}(\del_k x)^{\geq \ell +1} - \del_{k-1}(\del_k x)^{\leq \ell - 1}.
\]

\end{lem}
\begin{proof}
Since we have $\mathrm{Id}_{R_{k}(G)} = (-)^{\geq \ell+1} + (-)^{\leq \ell} = (-)^{\geq \ell+1} + (-)^{\ell} + (-)^{\leq \ell - 1}$ and $\del_{k-1}\circ \del_{k} = 0$, we obtain 
\begin{align*}
0 &= (\del_{k-1}\circ \del_{k})x \\
&= \del_{k-1}((\del_{k}x)^{\geq \ell+1} + (\del_{k}x)^{\ell} + (\del_{k}x)^{\leq \ell - 1}) \\
&= \del_{k-1}(\del_{k}x)^{\geq \ell+1} + \del_{k-1}(\del_{k}x)^{\ell} + \del_{k-1}(\del_{k}x)^{\leq \ell - 1},
\end{align*}
for every $x \in R_k(G)$. This completes the proof.
\end{proof}
\begin{df}\label{defmag1}
For $k, \ell \in \Z$, we define $\MC^{\ell}_{k}(G) = (R_k(G))^{\ell}$. That is, 
\[\MC^{\ell}_{k}(G) = \Z\langle (x_0, \dots, x_k) \in R_{k}(G) \mid L(x_0, \dots, x_k) = \ell\rangle \subset R_k(G),
\]
for $k \geq -1$, and  $\MC^{\ell}_{k}(G) = 0$ for $k \leq -2$.
We also define a homomorphism $\del^{\ell}_{k} : \MC^{\ell}_{k}(G) \too \MC^{\ell}_{k-1}(G)$ by $\del^{\ell}_{k}x = (\del_k x)^{\ell}$ for $x \in \MC^{\ell}_{k}(G)$.
\end{df}
\begin{lem}\label{magnitudehom}
$\del^{\ell}_{k-1}\circ \del^{\ell}_{k} = 0$.
\end{lem}
\begin{proof}
For every $x \in \MC^{\ell}_{k}(G)$, we have 
\begin{align*}
(\del^{\ell}_{k-1}\circ \del^{\ell}_{k})x &= (\del_{k-1}(\del_k x)^{\ell})^{\ell} \\
&= (-\del_{k-1}(\del_k x)^{\geq \ell +1} - \del_{k-1}(\del_k x)^{\leq \ell - 1})^{\ell} \\
& = 0.
\end{align*}
Here we applied Lemma \ref{easylemma} to the first line, and $(\del_k x)^{\geq \ell +1} = (\del_{k-1}(\del_k x)^{\leq \ell - 1})^{\ell} = 0$ to the second line. The latter can be verified by Lemma \ref{id}. This completes the proof.
\end{proof}
Lemma \ref{magnitudehom} implies that $(\MC^{\ell}_\ast(G), \partial^{\ell}_\ast)$ is a chain complex for each $\ell \in \Z$. 
\begin{df}\label{defmag2}
We call the homology of the chain complex $(\MC^{\ell}_\ast(G), \partial^{\ell}_\ast)$ \textit{magnitude homology} of a digraph $G$, and we denote it by $\MH^{\ell}_\ast(G)$.
\end{df}
\begin{rem}\label{embed}
For an undirected graph $g$, we can make it into a directed graph $\overrightarrow{g}$ by orienting each edge of $g$ in both directions. This operation is a functor from the category of undirected graphs and graph maps to that of digraphs and digraph maps, which is left adjoint to the forgetful functor. This forgetful functor sends a digraph to the undirected graph with the same vertex set, and where $v$ and $w$ share an edge if and only if the directed edges $(v, w)$ and $(w, v)$ are both present in the original digraph. We note that this left adjoint funtor embeds the category of undirected graphs into a full-subcategory of that of digraphs. The homology $\MH^{\ell}_k(\overrightarrow{g})$ defined above is isomorphic to the usual magnitude homology of an undirected graph $g$ defined in \cite{HW}.
\end{rem}
\begin{df}
For  $\ell, k \in \Z$, we define a homomorphism $\del'^{\ell}_{k} : \MC^{\ell}_{k}(G) \too \MC^{\ell-1}_{k-1}(G)$ by $\del'^{\ell}_{k}x = (\del_k x)^{\ell - 1}$.
\end{df}
\begin{lem}\label{kerker}
 $\del'^{\ell}_{k}(\mathrm{Ker}\ \del^{\ell}_{k}) \subset \mathrm{Ker}\ \del^{\ell-1}_{k-1}$.
\end{lem}
\begin{proof}
For every $x \in \mathrm{Ker}\ \del_k^{\ell}$, we have 
\begin{align*}
(\del^{\ell-1}_{k-1}\circ \del'^{\ell}_{k})x &= \del^{\ell-1}_{k-1}((\del_k x)^{\ell - 1}) \\
&= (\del_{k-1}((\del_k x)^{\ell - 1}))^{\ell-1} \\
&= (\del_{k-1}((\del_k x)^{\leq \ell - 1}))^{\ell-1} \\
&= (-\del_{k-1}(\del_k x)^{\ell} -\del_{k-1}(\del_k x)^{\geq \ell+1})^{\ell - 1} \\
&= 0.
\end{align*}
Here we applied $(\del_{k-1}(\del_kx)^{\leq \ell -2})^{\ell-1} = 0$ to the second line, Lemma \ref{easylemma} to the third line, and $(\del_k x)^{\ell} = (\del_k x)^{\geq \ell+1} = 0$ to the fourth line. These can be verified by Lemma \ref{id} or by the assumption.  This completes the proof.
\end{proof}
\begin{lem}\label{imim}
$\del'^{\ell}_k(\mathrm{Im}\ \del^{\ell}_{k+1}) \subset \mathrm{Im}\ \del^{\ell-1}_{k}$.
\end{lem}
\begin{proof}
For every $y \in \MC^{\ell}_{k+1}(G)$, we have 
\begin{align*}
\del'^{\ell}_{k}\del^{\ell}_{k+1}y &= \del'^{\ell}_{k}(\del_{k+1}y)^{\ell} \\
&= (\del_k(\del_{k+1}y)^{\ell})^{\ell - 1} \\
&= (-\del_k(\del_{k+1}y)^{\geq \ell+1} - \del_k(\del_{k+1}y)^{\leq \ell-1})^{\ell - 1} \\
&= (- \del_k(\del_{k+1}y)^{\ell-1})^{\ell - 1} \\
&= -\del^{\ell - 1}_{k}(\del_{k+1}y)^{\ell-1}.
\end{align*}
Here we applied Lemma \ref{easylemma} to the second line, and $(\del_{k+1}y)^{\geq \ell+1} = (\del_k(\del_{k+1}y)^{\leq \ell - 2})^{\ell-1} = 0$ to the third line. These can be verified by Lemma \ref{id}. This completes the proof.
\end{proof}
\begin{lem}\label{kerim}
$(\del'^{\ell - 1}_{k-1} \circ \del'^{\ell}_{k})(\mathrm{Ker}\ \del^{\ell}_{k}) \subset \mathrm{Im}\  \del^{\ell - 2}_{k-1}$.
\end{lem}
\begin{proof}
For every $x \in \mathrm{Ker}\ \del^{\ell}_{k}$, we have 
\begin{align*}
(\del'^{\ell - 1}_{k-1} \circ \del'^{\ell}_{k})x &= \del'^{\ell - 1}_{k-1}(\del_k x)^{\ell - 1} \\
&= (\del_{k-1}(\del_k x)^{\ell - 1})^{\ell - 2} \\
&= (\del_{k-1}((\del_k x)^{\leq \ell - 1} - (\del_k x)^{\leq \ell - 2}))^{\ell - 2} \\
&= (\del_{k-1}(-(\del_k x)^{\ell} -(\del_k x)^{\geq \ell+1} - (\del_k x)^{\ell - 2}))^{\ell - 2} \\
&= (- \del_{k-1}(\del_k x)^{\ell - 2})^{\ell - 2} \\
&= -\del^{\ell - 2}_{k-1} (\del_k x)^{\ell - 2}.
\end{align*}
Here we applied Lemma \ref{easylemma} and $(\del_{k-1}(\del_kx)^{\leq \ell - 3})^{\ell-2} = 0$ to the third line, and $(\del_k x)^{\ell} = (\del_k x)^{\geq \ell+1} = 0$ to the fourth line. These can be verified by Lemma \ref{id}. This completes the proof.
\end{proof}
Lemmas \ref{kerker} and  \ref{imim} implies that $\partial'^{\ell}_{k} : \MC^{\ell}_k(G) \too \MC^{\ell-1}_{k-1}(G)$ induces a homomorphism $\MH^{\ell}_{k}(G)\too \MH^{\ell-1}_{k-1}(G)$, which we also denote by $\del'^{\ell}_{k}$. Further, Lemma \ref{kerim} implies that $\del'^{\ell-1}_{k-1}\circ \del'^{\ell}_k = 0$, hence the pair $(\MH^{\ell - \ast}_{k - \ast}(G), \del'^{\ell - \ast}_{k - \ast})$ is a chain complex. 
\begin{df}\label{defdmh}
We denote the homology of  $(\MH^{\ell - \ast}_{k - \ast}(G), \del'^{\ell - \ast}_{k - \ast})$ by $\dMH^{\ell - \ast}_{k - \ast}(G)$.
\end{df}
\begin{eg}
It will turn out that $\dMH^k_k(G)$ is isomorphic to the $k$-th reduced path homology (Proposition \ref{pathdmh}). Hence $\dMH^0_0(G)$ is a free module, and its rank is equal to the number of connected components of $G$ (in the directed sense) minus 1. Proposition 2.9 of \cite{GLMY} shows that possible cycles of $\dMH^1_1(G)$ are made from double edges, triangles and squares in $G$.
\end{eg}
\section{Functoriality of magnitude homologies}\label{funct}
In this section, we show functoriality of $\MH^{\ell}_\ast(G)$ and $\dMH^{\ell}_\ast(G)$ after some definitions and technical lemmas.
\begin{df}
Let $G, H$ be digraphs. A \textit{digraph map} $f : G \too H$ is a map $f : V(G) \too V(H)$ satisfying that $(f(x), f(y)) \in E(H)$ or $f(x) = f(y)$ whenever $(x, y) \in E(G)$. Equivalently, a digraph map is a distance non-increasing map $V(G) \too V(H)$ with respect to the directed path metric.
\end{df}
\begin{df}
Let $f : G \too H$ be a digraph map and let $k \in \Z$. We define $(f_\#)_k : R_k(G) \too R_k(H)$ by the induced homomorphism of $f_\# : \Lambda_k(G) \too \Lambda_k(H)$ defined by  $(f_\#)_k (x_0, \dots, x_k) = (f(x_0), \dots, f(x_k))$. 
\end{df}
The following can be easily checked.
\begin{lem}\label{sharpchain}
$\{(f_\#)_k\}_k$ is a chain map.
\end{lem}
\begin{lem}\label{fdecomp}
For every $x \in R_k(G)$, we have $(f_{\#})_{k-1} \circ \del_k x = \del_k ((f_{\#})_{k} x)^{\geq \ell +1} + \del_k ((f_{\#})_{k}x)^{\ell} + \del_k ((f_{\#})_{k} x)^{\leq \ell -1})$.
\end{lem}
\begin{proof}
It is clear from Lemma \ref{sharpchain} and $\mathrm{Id}_{R_k(G)} = (-)^{\geq \ell+1} + (-)^{\ell} + (-)^{\leq \ell-1}$.
\end{proof}
\begin{df}
Let $f : G \too H$ be a digraph map. For $k, \ell \in \Z$, we define $(f_{\#})_{k}^{\ell} : R_k(G) \too R_k(H)$ by $(f_{\#})_{k}^{\ell} = ((f_{\#})_{k})^\ell$.
\end{df}
It is clear that $(f_{\#})_{k}^{\ell}$ is restricted to $\MC_{k}^{\ell}(G) \too \MC_{k}^{\ell}(H)$.
\begin{lem}\label{mccom}
For every $\ell \in \Z$, the family of homomorphisms $\{(f_{\#})_{k}^{\ell}\}$ restricted to $\MC_{k}^{\ell}(G) \too \MC_{k}^{\ell}(H)$ is a chain map. That is, the following diagram is commutative for every $k, \ell \in \Z \ :$
\begin{equation*}
    \xymatrix{
    \MC_{k}^{\ell}(G) \ar[r]^{(f_{\#})_{k}^{\ell}} \ar[d]_{\del_{k}^{\ell}} & \MC^{\ell}_{k}(H) \ar[d]^{\del_{k}^{\ell}}\\
    \MC_{k-1}^{\ell}(G) \ar[r]_{(f_{\#})_{k-1}^{\ell}} & \MC_{k-1}^{\ell}(H).
   }
\end{equation*}
\end{lem}
\begin{proof}
For every $x \in \MC^{\ell}_k(G)$, we have 
\begin{align*}
    \del_{k}^{\ell} \circ (f_{\#})_{k}^{\ell} x &= \del_{k}^{\ell} ((f_{\#})_{k} x)^{\ell} \\
    &= (\del_k ((f_{\#})_{k} x)^{\ell})^{\ell}\\
    &= ((f_{\#})_{k-1} \circ \del_k x - \del_k ((f_{\#})_{k} x)^{\geq \ell+1} - \del_k ((f_{\#})_{k} x)^{\leq \ell-1})^{\ell} \\
    &= ((f_{\#})_{k-1} \circ \del_k x)^{\ell} \\
    &= (f_\#)^{\ell}_{k-1}((\del_kx)^\ell + (\del_kx)^{\leq \ell - 1}) \\
    &= (f_{\#})_{k-1}^{\ell} \circ \del_{k}^{\ell}x.
\end{align*}
Here, we applied Lemma \ref{fdecomp} to the second line, $((f_\#)_kx)^{\geq \ell + 1} = (\del_k ((f_\#)_kx)^{\leq \ell - 1})^{\ell} = 0$ to the third line, and $(f_\#)^{\ell}_{k-1}(\del_kx)^{\leq \ell - 1} = 0$ to the fifth line. This completes the proof.
\end{proof}
\begin{lem}\label{mhcom}
For $x \in \mathrm{Ker}\ \del_{k}^{\ell} \subset \MC_{k}^{\ell}(G)$, we have $(f_{\#})_{k-1}^{\ell-1} \circ \del'^{\ell}_{k}x - \del'^{\ell}_{k}\circ (f_{\#})_{k}^{\ell}x \in \mathrm{Im}\ \del_{k}^{\ell-1}$.
\end{lem}
\begin{proof}
By the definition, we have $(f_{\#})_{k-1}^{\ell-1} \circ \del'^{\ell}_{k}x = ((f_{\#})_{k-1} (\del_k x)^{\ell - 1})^{\ell - 1}$ and $\del'^{\ell}_{k}\circ (f_{\#})_{k}^{\ell}x = (\del_k ((f_{\#})_{k} x)^{\ell })^{\ell - 1}$. Hence we have 
\begin{align*}
    (f_{\#})_{k-1}^{\ell-1} \circ \del'^{\ell}_{k}x - \del'^{\ell}_{k}\circ (f_{\#})_{k}^{\ell}x = \ &((f_{\#})_{k-1} (\del_k x)^{\ell - 1} - \del_k ((f_{\#})_{k} x)^{\ell })^{\ell - 1} \\
    = \ &((f_{\#})_{k-1} (\del_k x)^{\ell - 1} + \del_k ((f_{\#})_{k} x)^{\geq \ell + 1 } \\
    \ &+ \del_k ((f_{\#})_{k} x)^{\leq \ell - 1 } - (f_{\#})_{k-1} \circ \del_k x)^{\ell - 1} \\
    = \ &((f_{\#})_{k-1}( (\del_k x)^{\ell - 1} - \del_k x) +  \del_k ((f_{\#})_{k} x)^{\leq \ell - 1 } )^{\ell - 1}\\
    = \ &(-(f_{\#})_{k-1}( (\del_k x)^{\geq \ell + 1} + (\del_k x)^{\ell} + (\del_k x)^{\leq \ell-2}) \\
    &\ +  \del_k ((f_{\#})_{k} x)^{\leq \ell - 1 } )^{\ell - 1} \\
    = \ &(\del_k ((f_{\#})_{k} x)^{\leq \ell - 1 } )^{\ell - 1} \\
    = \ &\del_{k}^{\ell - 1}  ((f_{\#})_{k} x)^{\leq \ell - 1 } .
\end{align*}
Here we applied Lemma \ref{fdecomp} to the first line, Lemma \ref{easylemma} to the fourth line, and $(\del_kx)^{\geq \ell + 1} = (\del_kx)^{\ell } = ((f_\#)_{k-1}(\del_kx)^{\leq \ell - 2})^{\ell - 1} = 0$ to the fifth line. This completes the proof.
\end{proof}
\begin{df}
We define a homomorphism $(f_{\ast})_{k}^{\ell} : \MH_{k}^{\ell}(G) \too \MH_{k}^{\ell}(H)$ as the one induced from $(f_{\#})_{k}^{\ell}$.
\end{df}
\begin{prop}\label{dmhfunct}
The family of homomorphisms $\{(f_{\ast})_{k}^{\ell}\}$ is a chain map. That is , the following diagram is commutative for every $k, \ell \in \Z \ :$
\begin{equation*}
    \xymatrix{
    \MH_{k}^{\ell}(G) \ar[r]^{(f_{\ast})_{k}^{\ell}} \ar[d]_{\del'^{\ell}_{k}} & \MH^{\ell}_{k}(H) \ar[d]^{\del'^{\ell}_{k}}\\
    \MH_{k-1}^{\ell-1}(G) \ar[r]_{(f_{\ast})_{k-1}^{\ell-1}} & \MH_{k-1}^{\ell-1}(H).
    }
\end{equation*}
\end{prop}
\begin{proof}
It follows from Lemma \ref{mhcom}.
\end{proof}
\begin{df}
We define a homomorphism $(f^{1}_{\ast})_{k}^{\ell} : \dMH^{\ell}_{k}(G) \too \dMH^{\ell}_{k}(H)$ as the one induced from $(f_{\ast})_{k}^{\ell}$.
\end{df}
\section{An operator $\wedge$ and a product structure on complexes}\label{sectop}
In this section, we consider operations on chain complexes defined in the previous section. These operations are generalization of those introduced in \cite{GLMY}. The content of this section will be applied in Section \ref{sectinv}.
\begin{df}
For $k \in \Z$, we define a homomorphism $\bullet : \Lambda_i(G) \otimes \Lambda_j(G) \too \Lambda_{i+j+1}(G)$ by 
\[
\bullet ((x_0, \dots, x_i)\otimes (y_0, \dots, y_j)) = (x_0, \dots, x_i, y_0, \dots, y_j).
\]
We will use the notation $\bullet (u\otimes v) = u\bullet v = uv$ for simplicity.
\end{df}
\begin{lem}\label{prodcomm0}
The family of homomorphisms $\{\bullet : \Lambda_i(G) \otimes \Lambda_j(G) \too \Lambda_{i+j+1}(G)\}$ defines a degree $1$ chain map $\bullet : \Lambda(G)_\ast \otimes \Lambda(G)_\ast \too \Lambda_{\ast+1}(G)$.
\end{lem}
\begin{proof}
For $(x_0, \dots, x_i)\otimes (y_0, \dots, y_j) \in \Lambda_i(G)\otimes \Lambda_j(G)$, we have
\begin{align*}
    &\bullet \circ (\del \otimes \del)_{i+j}((x_0, \dots, x_i)\otimes (y_0, \dots, y_j)) \\
    = \ &\bullet(\del_i(x_0, \dots, x_i)\otimes (y_0, \dots, y_j) + (-1)^{i+1}(x_0, \dots, x_i)\otimes \del_j(y_0, \dots, y_j)) \\
    = \ &\del_i(x_0, \dots, x_i)\bullet (y_0, \dots, y_j) + (-1)^{i+1}(x_0, \dots, x_i)\bullet \del_j(y_0, \dots, y_j)
\end{align*}
We also have 
\begin{align*}
&\del_i(x_0, \dots, x_i)\bullet (y_0, \dots, y_j) \\ 
= \ &\sum_{p = 0}^{i-1}(-1)^p(x_0, \dots, \check{x}_p, \dots, x_i)\bullet (y_0, \dots, y_j) + (-1)^i(x_0, \dots, x_{i-1})\bullet (y_0, \dots, y_j) \\
= \ &\sum_{p = 0}^{i-1}(-1)^p(x_0, \dots, \check{x}_p, \dots, x_i, y_0, \dots, y_j) + (-1)^i(x_0, \dots, x_{i-1})\bullet (y_0, \dots, y_j),
\end{align*}
and 
\begin{align*}
&(-1)^{i+1}(x_0, \dots, x_i)\bullet \del_j(y_0, \dots, y_j) \\
= \ &(-1)^{i+1}(x_0, \dots, x_i)\bullet\sum_{q = 1}^{j}(-1)^q(y_0, \dots, \check{y}_q, \dots, y_j) + (-1)^{i+1}(x_0, \dots, x_{i})\bullet (y_1, \dots, y_j) \\
= \ & \sum_{q = 1}^{j}(-1)^{i+q+1}(x_0, \dots, x_i, y_0, \dots, \check{y}_q, \dots, y_j) + (-1)^{i+1}(x_0, \dots, x_{i})\bullet (y_1, \dots, y_j).
\end{align*}
Therefore we obtain
\begin{align*}
&\ \bullet \circ (\del \otimes \del)_{i+j}((x_0, \dots, x_i)\otimes (y_0, \dots, y_j)) \\
&= \del_i(x_0, \dots, x_i)\bullet (y_0, \dots, y_j) + (-1)^{i+1}(x_0, \dots, x_i)\bullet \del_j(y_0, \dots, y_j) \\
&=  \del_{i+j+1}(x_0, \dots, x_i, y_0, \dots, y_j)  \\
&=\del_{i+j+1}((x_0, \dots, x_i)\bullet (y_0, \dots, y_j)).
\end{align*}
This completes the proof.
\end{proof}
\begin{lem}\label{prodcomm}
The chain map $\bullet : \Lambda(G)_\ast \otimes \Lambda(G)_\ast \too \Lambda_{\ast+1}(G)$ induces a degree $1$ chain map $\bullet : R(G)_\ast \otimes R(G)_\ast \too R_{\ast+1}(G)$.
\end{lem}
\begin{proof}
Note that $\bullet$ is restricted the chain map $\Lambda_\ast(G) \otimes I_\ast(G) \oplus I_\ast(G) \otimes \Lambda_\ast(G) \too I_{\ast +1}(G)$. Hence the statement follows from Lemma \ref{prodcomm0} and the fact that $R_\ast(G) \otimes R_\ast(G)$ is isomorphic to the quotient of $ \Lambda_\ast(G)\otimes \Lambda_\ast(G)$ by $\Lambda_\ast(G) \otimes I_\ast(G) \oplus I_\ast(G) \otimes \Lambda_\ast(G)$. This completes the proof.
\end{proof}
Note that the homomorphism $\bullet : R(G)_i \otimes R(G)_j \too R_{i+j+1}(G)$ is explicitly defined as 
\[
(x_0, \dots, x_i)\bullet (y_0, \dots, y_j) = \begin{cases} (x_0, \dots, x_i, y_0, \dots, y_j) & x_i \neq y_0, \\ 0 & \text{otherwise}.\end{cases}
\]

\begin{df}
Let $G$ and $H$ be digraphs. \textit{The cartesian product} $G\times H$ is defined as a digraph with $V(G\times H) = V(G)\times V(H)$ and $((v, h), (v', h')) \in E(V\times H)$ if and only if $v= v'$ and $(h, h') \in E(H)$ or $(v, v') \in E(V)$ and $h = h'$.
\end{df}
Let $I$ be the digraph with $V(I) = \{0, 1 \}$ and $E(I) = \{(0, 1)\}$. We denote vertices $(v, 0)$ and $(v, 1)$ of $G\times I$ by $v$ and $v'$ respectively. Further, for $x = (x_0, \dots, x_k) \in R_k(G)$, we denote $(x'_0, \dots, x'_k) \in R_k(G\times I)$ by $x'$ and extend this notation linearly to $R_k(G)$. Note that this notation is well-defined as a homomorphism $\MC_k^{\ell}(G) \too \MC_k^{\ell}(G\times I)$, and induces a homomorphism $\MH_k^{\ell}(G) \too \MH_k^{\ell}(G\times I)$. The following operation is introduced in \cite{GLMY} as an analogy of the prism operator in singular homology theory.
\begin{df}
We define a homomorphism  $\wedge : R_k (G) \too R_{k+1} (G\times I)$ by 
\[
\wedge (x_0, \dots, x_k) = \sum_{i = 0}^{k} (-1)^i (x_0, \dots, x_i, x'_i, \dots, x'_k),
\]
for $k \geq 0$, and $\wedge = 0$ for $k \leq -1$. We will use the notation $\wedge u = \widehat{u}$.
\end{df}
Note that $\wedge$ is restricted to $ \MC^{\ell}_{k}(G) \too \MC^{\ell+1}_{k+1}(G\times I)$. 
The following two lemmas (\ref{hatprod} and \ref{delwedge}) appear in the proof of Proposition 2.12 in \cite{GLMY}.
\begin{lem}\label{hatprod}
Let $x \in R_p(G)$ and $y \in R_q(G)$. Then we have 
\[
\wh{xy} = \wh{x}y' + (-1)^{p+1}x\wh{y}.
\]

\end{lem}
\begin{proof}
It is easily checked for the case that $p \leq -1$ or $q \leq -1$. Hence it suffices to prove for the case that $x = (x_0, \dots, x_p)$ and $y = (y_0, \dots, y_q)$ with $x_p \neq y_0$. By definitions, we have 
\begin{align*}
\wh{xy} &= \sum_{i = 0}^{p}(-1)^i(x_0, \dots, x_i, x'_i, \dots, x'_p, y'_0, \dots, y'_q) + \sum_{j=0}^{q}(-1)^{p+1+j}(x_0, \dots, x_p, y_0, \dots, y_j, y'_j, \dots,  y'_q) \\
&= \wh{x}y' + (-1)^{p+1}x\wh{y}.
\end{align*}
This completes the proof.
\end{proof}
\begin{lem}\label{delwedge}
The homomorphism  $\wedge$ is a chain homotopy between chain maps $R_\ast(G) \too R_\ast(G\times I) ; x \mapsto x$ and $x \mapsto x'$. That is, for $x \in R_k(G)$, we have $\del_{k+1} \wh{x} = -\wh{\del_k x} + x' - x$.
\end{lem}
\begin{proof}
It is obvious for $k \leq -1$. Hence it suffices to prove for the case that $x = (x_0, \dots, x_k)$. We show that by induction on $k$. It is clear for $k = 0$. Suppose that we have $x = yx_{k+1}$ for some $y \in R_{k}(G)$ and $x_{k+1} \in V(G)$. Then we have
\begin{align*}
\del_{k+2} \wh{x} &= \del_{k+2}\wh{yx_{k+1}} \\
&= \del_{k+2}(\wh{y}x'_{k+1} + (-1)^{k+1}y(x_{k+1}, x'_{k+1})) \\
&= (\del_{k+1}\wh{y})x'_{k+1} + (-1)^{k+2}\wh{y} + (-1)^{k+1}(\del_k y)(x_{k+1}, x'_{k+1}) + y(x'_{k+1} - x_{k+1}) \\
&= (-\wh{\del_{k} y} + y' - y)x'_{k+1} + (-1)^{k+2}\wh{y} + (-1)^{k+1}(\del_k y)(x_{k+1}, x'_{k+1}) + y(x'_{k+1} - x_{k+1}) \\
&= -\wh{\del_{k} y}x'_{k+1} +  (-1)^{k+1}(\del_k y)(x_{k+1}, x'_{k+1}) +  (-1)^{k+2}\wh{y} + x'  - x  \\
&= -\wh{(\del_k y)x_{k+1}} + (-1)^{k+2}\wh{y} + x' - x \\
&= -\wedge((\del_k y)x_{k+1} + (-1)^{k+1}y) + x' - x \\
&= -\wh{\del_{k+1}x} + x' - x.
\end{align*}
Here we applied Lemma \ref{hatprod} to the first line, Lemma \ref{prodcomm} to the second line, and the induction assumption to the third line. This completes the proof.
\end{proof}
\begin{lem}\label{lwedge}
For every $k, \ell \in \Z$ and $x \in R_k(G)$, we have $(\wh{x})^{ \ell + 1} = \wh{x^{\ell}}$.
\end{lem}
\begin{proof}
It is sufficient to consider the case that $x = (x_0, \dots, x_k)$ by the linearity of $\wedge$ and $(-)^{ \ell}$. Note that $L(\wh{x}) = L(x) + 1$ and $x^\ell = \begin{cases} x & L(x) = \ell, \\ 0 & \text{otherwise}\end{cases}$. Now we have 
\[
(\wh{x})^{\ell + 1}  = \begin{cases} \wh{x} & L(x) = \ell, \\ 0 & \text{otherwise},\end{cases} = \wedge \begin{cases} x & L(x) = \ell, \\ 0 & \text{otherwise},\end{cases} = \wh{x^\ell}.
\]
 This completes the proof.
\end{proof}
\begin{prop}\label{wedgechain}
Let us consider a homomorphism $(-1)^k \wedge : \MC_{k}^{\ell}(G) \too \MC_{k+1}^{\ell+1}(G\times I)$. The family of homomorphisms $\{(-1)^k \wedge\}$ defines a degree $1$ chain map $\MC_{k + \ast}^{\ell + \ast}(G) \too \MC_{k + \ast+1}^{\ell + \ast+1}(G\times I)$.
\end{prop}
\begin{proof}
For $x \in \MC_{k}^{\ell}(G)$, we have 
\begin{align*}
    \del_{k+1}^{\ell + 1}(-1)^k \wh{x} &= (-1)^k (\del_{k+1}\wh{x})^{\ell + 1} \\
    &= (-1)^k (-\wh{\del_{k}x} + x' - x)^{\ell + 1} \\
    &= (-1)^{k-1} (\wh{\del_{k}x})^{\ell + 1} \\
    &= (-1)^{k-1} \wh{(\del_{k}x)^{\ell}} \\
    &= (-1)^{k-1} \wh{\del_{k}^{\ell}x}.
\end{align*}
Here we applied Lemma \ref{delwedge} to the second line, and Lemma \ref{lwedge} to the third line. This completes the proof.
\end{proof}
\begin{cor}\label{wedgemh}
The homomorphism $\wedge : \MC_{k}^{\ell}(G) \too \MC_{k+1}^{\ell+1}(G\times I)$ induces a homomorphism $\wedge : \MH_{k}^{\ell}(G) \too \MH_{k+1}^{\ell+1}(G\times I)$. Further, this is a chain homotopy between chain maps $\MH_\ast^{\ast}(G) \too \MH_\ast^{\ast}(G\times I) ; x \mapsto x$ and $x \mapsto x'$.
\end{cor}
\begin{proof}
The former is clear from Proposition \ref{wedgechain}. The latter follows from the following : for $x \in \MC_k^{\ell}(G)$, we have 
\begin{align*}
    \del'^{\ell + 1}_{k+1}\wh{x} &= (\del_{k+1}\wh{x})^\ell \\
    &= -(\wh{\del_k x})^\ell + (x' - x)^{\ell} \\
    &= -\wh{(\del_k x)^{\ell -1}} + x' - x \\
    &= -\wh{\del'^{\ell}_k x} + x' - x.
\end{align*}
Here we applied Lemma \ref{delwedge} to the first line, and Lemma \ref{lwedge} to the second line. Hence the statement follows by taking homology. This completes the proof.
\end{proof}

\section{Homotopy invariance of $\dMH^{\ell}_{k}(G)$}\label{sectinv}
In this section, we show homotopy invariance of $\dMH^{\ell}_k(G)$, which will turn out to be a generalization of homotopy invariance of path homology shown in \cite{GLMY}. The following definition of the homotopy is introduced in \cite{GLMY}.
\begin{df}[\cite{GLMY} Section 3.1 ]
A digraph $G$ is a  \textit{line graph} of length $n$ if $V(G) = \{0, \dots, n\}$, $\#E(G) = n$, and exactly one of $(i, i+1)$ or $(i+1, i)$ is in $E(G)$ for each $0 \leq i \leq n-1$. We denote the set of all line graphs of length $n$ by $\mathcal{I}_n$. 
\end{df}
Note that we have $I  \in \mathcal{I}_1$ for the digraph $I$ defined in Section \ref{sectop}.
\begin{df}[\cite{GLMY} Definition 3.1]\label{defhtpy}
Let $f, g : G \too G'$ be digraph maps. A \textit{homotopy} between $f$ and $g$ is a digraph map $F : G\times I_n \too G'$ for some $n \geq 1$ and $I_n \in \mathcal{I}_n$ such that $F|_{G\times\{0\}} = f$ and $F|_{G\times\{n\}} = g$. If there exists a homotopy between digraph maps $f$ and $g$, we say they are \textit{homotopic} and denote it by $f \simeq g$.
\end{df}
It is easy to check that the relation $\simeq$ above  is an equivalence relation.
\begin{rem}[\cite{GLMY} Section 3.1]\label{onestep}
For $x, y \in G$, we use the notation $x \overrightarrow{=} y$ if $(x, y) \in E(G)$ or $x = y$. It is easily checked by the definition that $f, g : G \too G'$ are one-step homotopic, namely there is a homotopy $F : G\times I_1 \too G'$ between them, if and only if $f(x) \overrightarrow{=} g(x)$ for every $x \in V(G)$ or $g(x) \overrightarrow{=} f(x)$ for every $x \in V(G)$.
\end{rem}
\begin{df}
Digpraphs $G$ and $G'$ are \textit{homotopy equivalent} if there exist digraph maps $f : G \too G'$ and $g : G' \too G$ satisfying $g\circ f \simeq \mathrm{Id}_{G}$ and $f\circ g \simeq \mathrm{Id}_{G'}$. We say a digraph $G$ is \textit{contractible} if it is homotopy equivalent to the one point graph $\{0\} \in \mathcal{I}_0$.
\end{df}
\begin{eg}[\cite{GLMY} Example 3.10]\label{egcontr2}
 A digraph $G$ is a \textit{directed tree} if it has no undirected cycles. A digraph obtained from a directed tree $G$ by removing a leaf edge is homotopy equivalent to $G$ itself by \textit{Remark} \ref{onestep}. Hence any directed trees are contractible.
\end{eg}
\begin{eg}[\cite{GLMY} Example 3.11]\label{egcontr}
Let $K_n$ be the complete undirected graph with $n$ vertices. Then the identity map on $\overrightarrow{K_n}$ is homotopic to a one point map by \textit{Remark} \ref{onestep}. Hence  $\overrightarrow{K_n}$ is contractible. In general, a \textit{star-like} ( or a \textit{inverse star-like}) digraph is contractible, where a digraph $G$ is star-like ( or inverse star-like) if there exists a vertex $a \in V(G)$ such that $(a, x) \in E(G)$ for every $a \neq x \in V(G)$ ( or $(x, a) \in E(G)$ respectively). 
\end{eg}
\begin{prop}\label{htpyinv}
Let $f, g : G \too G'$ be digraph maps. If they are homotopic, then the induced homomorphisms $(f^{1}_\ast)_{k}^{\ell}, (g^{1}_\ast)^{\ell}_{k} : \dMH^{\ell}_{k}(G) \too \dMH^{\ell}_{k}(G')$ are identical, that is, we have $(f^{1}_\ast)_{k}^{\ell} = (g^{1}_\ast)_{k}^{\ell}$ for all $\ell, k \in \Z$. 
\end{prop}
For the proof, we use the following fundamental lemma of homological algebra, which we omit the proof.
\begin{lem}\label{chainhtpy}
Let $G : C_\ast \too D_{\ast + 1}$ be a chain homotopy between chain maps $f, g : C_\ast \too D_{\ast}$. For any chain map $F : C_\ast \too D_{\ast}$, $F\circ G$ is a chain homotopy between $F\circ f$ and $F\circ g$.
\end{lem}
\begin{proof}[Proof of Proposition \ref{htpyinv}]
Let $F : G\times I \too G'$ be a homotopy between $f$ and $g$. We construct a chain homotopy $L^{\ell + \ast}_{k + \ast} : \MH^{\ell + \ast}_{k + \ast}(G) \too \MH^{\ell + \ast + 1}_{k + \ast + 1}(G')$, that is a homomorphism satisfying $L^{\ell-1}_{k-1}\circ \del'^{\ell}_{k} + \del'^{\ell + 1}_{k+1}\circ L^{\ell}_{k} = (g_\ast)^{\ell}_{k} - (f_\ast)^{\ell}_{k}$ for all $\ell, k \in \Z$. It is an elementary fact that an existence of such a chain homotopy guarantees that $(f^{1}_\ast)_{k}^{\ell} = (g^{1}_\ast)_{k}^{\ell}$. We define $L^{\ell}_k : \MH^{\ell}_{k}(G) \too \MH^{\ell+1}_{k+1}(G')$ by $L^{\ell}_k x = (F_\ast)_{k+1}^{\ell+1} \wh{x}$. Then it follows from Cororally \ref{wedgemh} and Lemma \ref{chainhtpy} that $L^{\ast}_\ast$ is a desired chain homotopy. This completes the proof.


\end{proof}
\begin{rem}
 In Section \ref{sectcoin}, it will turn out that Proposition \ref{htpyinv} is a generalization of Theorem 3.3 in \cite{GLMY}.
\end{rem}
The following is clear from Proposition \ref{htpyinv}.
\begin{prop}\label{htpyequiv}
If two digraphs $G, G'$ are homotopy equivalent, then we have $\dMH^{\ell}_k(G) \cong \dMH^{\ell}_k(G')$ for every $\ell, k \in \Z$.
\end{prop}
\section{Coincidence with the path homology}\label{sectcoin}
In this section, we introduce Grigor'yan--Lin--Muranov--S-T.Yau's \textit{reduced path homology $\tilde{H}_\ast (G)$ of a digraph $G$} (\cite{GLMY}), and we show that $ \tilde{H}_k (G) \cong \dMH_{k}^{k}(G)$.

\begin{df}
Let $k \in \Z$. We define a module $A_{k}(G)$ by 
\[
A_k(G) = \Z\langle (x_0, \dots, x_k) \in V(G)^{k+1} \mid (x_i, x_{i+1}) \in E(G) \ \text{for}\  0\leq i\leq k-1\rangle,
\]
for $k \geq -1$ and $A_{k}(G) = 0$ for $k\leq -2$.
\end{df}
An element of $A_k(G)$ is called \textit{an allowed $k$-path} in \cite{GLMY}. The following is obvious.
\begin{lem}\label{amc}
$A_k(G) = \MC_{k}^{k}(G)$.
\end{lem}
By Lemma \ref{amc}, we can consider $A_k(G)$ as a submodule of $R_k(G)$.
\begin{df}\label{delreg}
Let $k \in \Z$. We define a homomorphism $\del^{\rm reg}_{k} : A_{k}(G) \too R_{k-1}(G)$ by $\del^{\rm reg}_{k} = \del_k|_{A_k(G)}$. 
\end{df}

\begin{df}
For $k \in \Z$, we define a module $\Omega_{k}^{\rm reg}(G)$ by 
\[
\Omega_{k}^{\rm reg}(G) = \{x \in A_k(G) \mid \del^{\rm reg}_{k}x \in A_{k-1}(G)\}.
\]
\end{df}
Note that we have $\Omega_{k}^{\rm reg}(G) \subset A_k(G) = \MC_{k}^{k}(G)$.
The following is clear from Lemma \ref{rchain}.
\begin{prop}
$(\Omega_{\ast}^{\rm reg}(G), \del^{\rm reg}_{\ast})$ is a chain complex.
\end{prop}
\begin{eg}\label{egkn}
For the undirected complete graph $K_n$ and any digraph $G$ with $\#V(G) = n$, we have $\Omega_{k}^{\rm reg}(\overrightarrow{K_n}) = A_k(\overrightarrow{K_n}) =  R_k(G)$. Hence we have $(\Omega_{\ast}^{\rm reg}(\overrightarrow{K_n}), \del^{\rm reg}_{\ast}) = (R_{\ast}(G), \del_\ast)$ for any digraph $G$ with $\#V(G) = n$.
\end{eg}
\begin{df}
We define $\tilde{H}_k(G) := H_k(\Omega_{\ast}^{\rm reg}(G), \del^{\rm reg}_{\ast})$ 
\end{df}
\begin{lem}\label{omeker}
$\Omega_{k}^{\rm reg}(G) = \mathrm{Ker}\ \del_{k}^{k} = \MH_{k}^{k}(G)$.
\end{lem}
\begin{proof}
Let $x \in A_k(G) = \MC_{k}^{k}(G)$. Then $x \in \Omega_{k}^{\rm reg}(G)$ if and only if we have $\del^{\rm reg}_{k}x = (\del_{k}x)^{k} + (\del_{k}x)^{k-1} \in A_{k-1}(G)$, which is equivalent to $\del^{k}_{k}x = (\del_{k}x)^{k} = 0$. This completes the proof.
\end{proof}
\begin{lem}\label{kerker2}
$\mathrm{Ker}\ \del^{\rm reg}_{k} = \mathrm{Ker}\ \del'^{k}_{k} \subset \MH_{k}^{k}(G).$
\end{lem}
\begin{proof}
Let $x \in \Omega^{\rm reg}_{k}(G) = \mathrm{Ker}\ \del^{k}_{k}$. Then  $x \in \mathrm{Ker}\ \del^{\rm reg}_{k}$ if and only if $\del^{\rm reg}_{k}x = (\del_{k}x)^{k-1} = 0$, which is equivalent to $x \in \mathrm{Ker}\ \del'^{k}_{k}$. This completes the proof.
\end{proof}
\begin{lem}
$ \del^{\rm reg}_{k+1} \Omega^{\rm reg}_{k+1}(G) = \del'^{k+1}_{k+1}\MH^{k}_{k}(G)$.
\end{lem}
\begin{proof}
Let $x \in \Omega^{\rm reg}_{k+1}(G) = \MH^{k+1}_{k+1}(G)$. Then we have $ \del^{\rm reg}_{k+1}x =  (\del_{k+1}x)^{k} = \del'^{k+1}_{k+1}x$. This completes the proof.
\end{proof}
Therefore, we obtain the following.
\begin{prop}\label{pathdmh}
We have 
\[
\tilde{H}_k(G) = \mathrm{Ker}\ \del^{\rm reg}_{k}/\del^{\rm reg}_{k+1} \Omega^{\rm reg}_{k+1}(G) \cong \mathrm{Ker}\ \del'^{k}_{k}/\del'^{k+1}_{k+1}\MH^{k}_{k}(G) = \dMH^{k}_{k}(G).
\]

\end{prop}
According to Proposition \ref{pathdmh}, we can consider homologies $\dMH^{\ell}_k(G)$ as a generalization of the path homology. As an application, we have the following.
\begin{prop}\label{rknull}
The homology of the chain complex $(R_\ast(G), \del_\ast)$ is all $0$ for any digraph $G$.
\end{prop}
\begin{proof}
By Example \ref{egkn}, we have $H_k(R_\ast(G)) =  \tilde{H}_k(\overrightarrow{K_n})$ for $n = \#V(G)$. Further, we have $\tilde{H}_k(\overrightarrow{K_n}) \cong \dMH^{k}_k(\overrightarrow{K_n})$ by Proposition \ref{pathdmh}. Since the digraph  $\overrightarrow{K_n}$ is contractible by Example \ref{egcontr}, we obtain that $\dMH^{k}_k(\overrightarrow{K_n}) \cong \dMH^{k}_k(\{0\}) \cong 0$ for all $k \in \Z$ by Proposition \ref{htpyequiv}. This completes the proof.
\end{proof}
\section{On a spectral sequence producing magnitude homologies}\label{sectspec}
In this section, we construct a spectral sequence whose first page is identical to magnitude homology. Further, the `diagonal part' of the second page is identical to path homology. 
\subsection{Preliminary for spectral sequences}
We first give a short introduction to spectral sequences. For details, see \cite{Wei} section 5.4. Readers familiar to this subject can skip this subsection, however, we slightly change notations from traditional one to make the correspondence with magnitude homologies clear, for example `$E_{p, q}^r$' to $E^{p, r}_{p+q}$. 

Let $(C_{\ast}, \del_{\ast})$ be a chain complex, that is a family of abelian groups $\{C_n\}_{n \in \Z}$ with homomorphisms $\del_n : C_n \too C_{n-1}$ satisfying $\del_{n+1}\circ \del_n = 0$. We suppose that $(C_{\ast}, \del_{\ast})$ is filtered, that is, we have 
\[
C_n = \bigcup_{\ell \in \Z} C^{\ell}_{n}
\]
with 
\[
\cdots \subset C^{\ell}_n \subset C^{\ell + 1}_{n} \subset \cdots \subset C_n
\]
for each $n \in \Z$, and $\del_{n+1}C^{\ell}_{n+1} \subset C^{\ell}_n$ for each $n, \ell \in \Z$. We further suppose that this filtration is bounded, that is, for each $n \in \Z$ there exists $L >0$ with $C^{L}_{n} = C_n$ and $C^{-L}_{n} = 0$. We define the following submodules of $C^{\ell}_n$ : 
\begin{align*}
    Z^{\ell, r}_n &:= \{x \in C^{\ell}_n \mid \del_nx \in C^{\ell-r}_{n-1}\}, \\
    B^{\ell, r}_n &:= \{x \in C^{\ell}_{n} \mid \exists y \in C^{\ell+r}_{n+1} \ \del_{n+1}y = x\}, \\
    Z^{\ell, \infty}_n &:= \{x \in C^{\ell}_n \mid \del_nx = 0\}, \\
    B^{\ell, \infty}_n &:= \{x \in C^{\ell}_{n} \mid \exists y \in C_{n+1}\  \del_{n+1}y = x\}.
\end{align*}
Note that we have 
\[
B^{\ell, 0}_n \subset \cdots \subset B^{\ell, r}_n \subset B^{\ell, r+1}_n \subset \cdots \subset B^{\ell, \infty}_n \subset Z^{\ell, \infty}_n \subset \cdots \subset Z^{\ell, r+1}_n \subset Z^{\ell, r}_n \subset \cdots \subset Z^{\ell, 0}_n = C^{\ell}_n,
\]
and $Z^{\ell-1, r-1}_n \subset Z^{\ell, r}_n$. Hence we have $ Z^{\ell-1, r-1}_n + B^{\ell, r-1}_n \subset  Z^{\ell, r}_n$, and we define 
\[
E^{\ell, r}_n := Z^{\ell, r}_n/Z^{\ell-1, r-1}_n + B^{\ell, r-1}_n,
\]
and 
\[
E^{\ell, \infty}_n := Z^{\ell, \infty}_n/Z^{\ell-1, \infty}_n + B^{\ell, \infty}_n.
\]
The following lemma is straightforward.
\begin{lem}\label{sp1}
$\del_n Z^{\ell+r, r}_n = B^{\ell, r}_{n-1}$
\end{lem}
The above Lemma \ref{sp1} implies that $\del_n$ induces a homomorphism 
\[
\del^{\ell, r}_n : E^{\ell, r}_n \too E^{\ell-r, r}_{n-1}
\]
with $\del^{\ell+r, r}_{n+1} \circ \del^{\ell, r}_n = 0$. The following lemmas are also straightforward.
\begin{lem}\label{sp2}
${\rm Ker}\ \del^{\ell, r}_n = Z^{\ell-1, r-1}_n + Z^{\ell, r+1}_n / Z^{\ell-1, r-1}_n + B^{\ell, r-1}_n$.
\end{lem}
\begin{lem}\label{sp3}
$\del^{\ell+r, r}_{n+1}E^{\ell+r, r}_{n+1} = B^{\ell, r}_n/Z^{\ell-1, r-1}_n + B^{\ell, r-1}_n$.
\end{lem}
The above Lemmas \ref{sp2} and \ref{sp3} implies that 
\[
{\rm Ker}\ \del^{\ell, r}_n/\del^{\ell+r, r}_{n+1}E^{\ell+r, r}_{n+1} \cong E^{\ell, r+1}_{n}.
\]
The next lemma follows from the assumption of boundedness.
\begin{lem}\label{spconv}
For each $\ell, n \in \Z$, there exits $L>0$ such that $E^{\ell, r}_{n} \cong E^{\ell, \infty}_n$ for every $r > L$.
\end{lem}
Note that the filtration on $(C_\ast, \del_{\ast})$ induces a filtration on the homology $H_n(C_\ast)$ as 
\[
F_\ell H_n(C_\ast) = \{x \in H_n(C_\ast) \mid \exists x' \in C^{\ell}_n\ [x'] = x\}.
\]
We define that $G_\ell H_n(C_\ast) := F_\ell H_n(C_\ast)/F_{\ell-1} H_n(C_\ast)$. The following is immediate.
\begin{lem}\label{spgr}
$E^{\ell, \infty}_n \cong G_\ell H_n(C_\ast).$
\end{lem}
\begin{rem}Note that in a traditional convention, we have `$E_{p, q}^r$' $= E^{p+q, r}_p$ and `$E_{p, q}^{\infty}$' $= E^{p+q, \infty}_p$.
\end{rem}
Next we consider the functoriality of spectral sequences. Let $C_\ast$ and $D_\ast$ be filetered chain complexes. We denote the corresponding modules considered above by $E^{\ell, r}_{n}(C_\ast)$ and $E^{\ell, r}_{n}(D_\ast)$ respectively. Let  $f : C_\ast \too D_\ast$ be a filtered chain map, that is a family of homomorphisms $f_n : C_n \too D_n$ satisfying $f_n\circ \del_n = \del_n \circ f_n$ and $f_n C^{\ell}_n \subset D^{\ell}_n$. Then it is easily checked that $f$ induces homomorphisms $f^{\ell, r}_n : E^{\ell, r}_{n}(C_\ast) \too E^{\ell, r}_{n}(D_\ast)$ which commutes with the differentials $\del^{\ell, r}_n$, namely a chain map. Further, Lemmas \ref{sp2} and \ref{sp3} implies that $f^{\ell, r}_n$ induces $f^{\ell, r+1}_n$ by taking homologies of $E^{\ell, r}_{n}(C_\ast)$ and $E^{\ell, r}_{n}(D_\ast)$. Then the following is obvious.
\begin{lem}\label{spfunc}
If $f^{\ell, r}_n$ is an isomorphism for all $n$, then so are $f^{\ell, r+k}_n$ for $k>0$.
\end{lem}

\subsection{Construction of the spectral sequence}

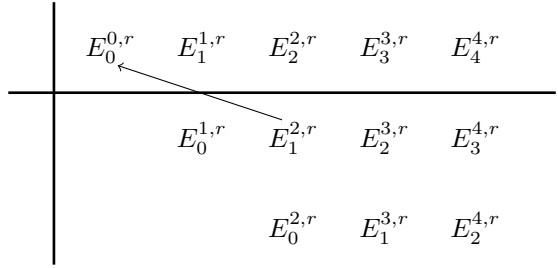
\begin{figure}[htbp]
\centering
\begin{tikzpicture}[scale=1.2]
\filldraw[fill=black, draw=black] (0, 0)  node[left] {$E^{0,r}_0$} ;
\filldraw[fill=black, draw=black] (1, 0)  node[left] {$E^{1,r}_1$} ;
\filldraw[fill=black, draw=black] (2, 0)  node[left] {$E^{2,r}_2$} ;
\filldraw[fill=black, draw=black] (3, 0)  node[left] {$E^{3,r}_3$} ;
\filldraw[fill=black, draw=black] (4, 0)  node[left] {$E^{4,r}_4$} ;

\filldraw[fill=black, draw=black] (1, -1)  node[left] {$E^{1,r}_0$} ;
\filldraw[fill=black, draw=black] (2, -1)  node[left] {$E^{2,r}_1$} ;
\filldraw[fill=black, draw=black] (3, -1)  node[left] {$E^{3,r}_2$} ;
\filldraw[fill=black, draw=black] (4, -1)  node[left] {$E^{4,r}_3$} ;

\filldraw[fill=black, draw=black] (2, -2)  node[left] {$E^{2,r}_0$} ;
\filldraw[fill=black, draw=black] (3, -2)  node[left] {$E^{3,r}_1$} ;
\filldraw[fill=black, draw=black] (4, -2)  node[left] {$E^{4,r}_2$} ;

\draw[line width=1pt] (-1,0.5)--(-1,-2.4);
\draw[line width=1pt] (-1.5,-0.5)--(4.5,-0.5);
\draw[->] (1.5,-0.8)--(-0.3,-0.2);
\end{tikzpicture}
\label{graphty}
\caption{Illustration of the spectral sequence $E^{\ell, r}_n$ in the $(\ell, n-\ell)$ coordinate , following the ordinary coordinate of `$E^r_{\ell, n-\ell}$'.}
\label{spe1}
\end{figure}

Now we consider the chain complex $(R_\ast(G), \del_\ast)$ for some digraph $G$ with a filtration 
\begin{align}\label{filt}
    R^{\ell}_{k}(G) := (R_k(G))^{\leq \ell},
\end{align}
that is, $R^{\ell}_k(G) = \Z\langle(x_0, \dots, x_k) \in V(G)^{k+1} \mid L(x_0, \dots, x_k) \leq \ell  , x_i \neq x_{i+1}\rangle$ for $k\geq -1$ and $R^{\ell}_k(G) = 0$ for $k\leq -2$. 
\begin{prop}\label{magspec}
The filtration (\ref{filt}) induces a spectral sequence such that 
\begin{align*}
    E^{\ell, 0}_k &= \MC^{\ell}_{k}(G), \\
    E^{\ell, 1}_k &= \MH^{\ell}_k(G), \\
    E^{\ell, 2}_k &=  \dMH^{\ell}_{k}(G).
\end{align*}
\end{prop}
\begin{proof}
We have 
\[
E^{\ell, 0}_k = R^{\ell}_k(G)/R^{\ell-1}_k(G) \cong \MC^{\ell}_{k}(G),
\]
and $\del^{\ell, 0}_{k} : E^{\ell, 0}_k \too E^{\ell, 0}_{k-1}$ coincides with $\del^{\ell}_k : \MC^{\ell}_{k}(G) \too \MC^{\ell}_{k-1}(G)$. Hence we obtain $E^{\ell, 1}_k = \MH^{\ell}_k(G)$. Now we see that $\del^{\ell, 1}_{k} : E^{\ell, 1}_k \too E^{\ell-1, 1}_{k-1}$ coincides with $\del'^{\ell}_k : \MH^{\ell}_{k}(G) \too \MH^{\ell-1}_{k-1}(G)$. It follows from the fact that both of $\del^{\ell, 1}_{k}x$ and $\del'^{\ell}_{k}x$ are represented by $\del_kx$ for $x \in E^{\ell, 1}_k =  \MH^{\ell}_{k}(G)$. Thus we obtain $E^{\ell, 2}_k =  \dMH^{\ell}_{k}(G)$. 
\end{proof}
\begin{prop}\label{diambd}
 For a digraph $G$ with $d := {\rm diam}\  G < \infty$, we have $E^{\ell, r}_n = 0$ for $n, \ell \in \Z$ with $nd < \ell$ and any $r \in \Z$.
\end{prop}
\begin{proof}
 It follows from that $E^{\ell, 0}_n = R^\ell_n(G)/R^{\ell-1}_n(G) = 0$ for such $n, \ell, r$'s.
\end{proof}
In this paper, we illustrate the above spectral sequence in the `ordinary' coordinate as Figure \ref{spe1}.
\section{Applications}\label{sectapl}
In this section, we give applications of our construction of $\dMH^{\ell}_k(G)$ and the spectral sequence.
\subsection{Homotopy invariance of $E^{\ell, r}_k$ for $r \geq 2$}
Proposition \ref{htpyequiv} and Proposition \ref{magspec} implies that $E^{\ell, 2}_k(G)$ is homotopy invariant for every $k \in \Z$. Together with functoriality of spectral sequences, we have the following. For digraphs $G, G'$, we denote the modules $E^{\ell, r}_k$ obtained from chain complexes $R_\ast(G)$ and $R_\ast(G')$ with the filtration (\ref{filt}) by $E^{\ell, r}_k(G)$ and $E^{\ell, r}_k(G')$ respectively.
\begin{prop}
 $E^{\ell, r}_k(G)$ is homotopy invariant for $r \geq 2$. That is, a homotopy equivalence $f : G \too G'$ induces isomorphisms $E^{\ell, r}_k(G) \cong E^{\ell, r}_k(G')$ for every $\ell, k \in \Z$ and $r \geq 2$.
\end{prop}
\begin{proof}
Let $f : G \too G'$ be a digraph map.  By Lemma \ref{sharpchain}, $f$ induces a chain map $f_\# : R_\ast(G) \too R_\ast(H)$. It is easily checked that $f_\#$ respects the filtration (\ref{filt}). Hence $f$ induces a homomorphism $f^{\ell, r}_{k} : E^{\ell, r}_k(G) \too E^{\ell, r}_k(H)$. Now we suppose that $f$ is a homotopy equivalence map. Then Lemma \ref{htpyinv} implies that $f^{\ell, 2}_k$ is an isomorphism for any $\ell, k$. Hence Lemma \ref{spfunc} implies the statement.
\end{proof}
\subsection{Diagonality implies vanishing of path homologies}
Now we suppose that a digraph $G$ has finite diameter, that is, we have $d(x, y) < \infty$ for any $x, y \in V(G)$. This assumption implies the boundedness of the filtration (\ref{filt}). Note that $\overrightarrow{g}$ satisfies this assumption for any undirected connected graph $g$.
\begin{lem}\label{zeroconv}
 Under the above finiteness condition, the spectral sequence induced from the filtration converges to $0$. That is, for every $\ell, k \in \Z$, $E^{\ell, r}_k(G) \cong E^{\ell, \infty}_k(G) \cong 0$ with sufficiently large $r$.
\end{lem}
\begin{proof}
 By Lemmas \ref{spconv} and \ref{spgr}, we have $E^{\ell, r}_k(G) \cong E^{\ell, \infty}_k(G) \cong G_\ell H_k(R_\ast)$ for sufficiently large $r$. Further, we have $H_k(R_\ast) \cong 0$ for every $k \in \Z$ by Proposition \ref{rknull}, hence we obtain $G_\ell H_k(R_\ast) \cong 0$. This completes the proof.
\end{proof}
The following notion of diagonality is introduced in \cite{HW} for undirected graphs.
\begin{df}
 A digraph $G$ is \textit{diagonal} if $\MH^{\ell}_k(G) = 0$ for $\ell \neq k$.
\end{df}
\begin{eg}\label{diageg}
\begin{enumerate}
\item  For any undirected tree $T$ and any undirected complete graph $K_n$,  $\overrightarrow{T}$ and $\overrightarrow{K_n}$ are diagonal (\cite{HW} Example 2.5 and Corollary 6.8). 
\item For undirected graphs $G$ and $H$, $G\star H$ denotes their join product. Then $\overrightarrow{G\star H}$ is diagonal for any $G$ and $H$ (\cite{HW} Theorem 7.5).
\item Let $G$ be the {\it icosahedral graph}. Then $\overrightarrow{G}$ is diagonal (\cite{Gu} Theorem 4.5).
\item A graph $G$ is called {\it pawful} if its diameter is $\leq 2$ and for any vertices $x, y, z$ with $d(x,y) = 1, d(x, z) = d(y, z) = 2$, there exists a vertex $w$ with $d(x, w) = d(y, w) = d(z, w) = 1$ (\cite{Gu} Definition 4.2). For any pawful graph $G$, $\overrightarrow{G}$ is diagonal (\cite{Gu} Theorem 4.4, \cite{TY} Corollary 3.5).
\item For undirected graphs $G$ and $H$, $G\times H$ denotes their cartesian product. Then $\overrightarrow{G\times H}$ is diagonal for any diagonal graphs $\overrightarrow{G}$ and $\overrightarrow{H}$ (\cite{HW} Proposition 7.3).
\end{enumerate}
\end{eg}
\begin{rem}
 The diagonality of other graphs is studied, for example, in \cite{TY}. It is known that the girth of any diagonal graph which is not a tree is $3$ or $4$ (\cite{AHK} Corollary 1.6).
\end{rem}
The following proposition relates magnitude homology and path homology.
\begin{prop}\label{diagvanish}
 If a digraph $G$ is diagonal and has finite diameter, then the reduced path homology $\tilde{H}_k(G)$ is $0$ for every $k \in \Z$.
\end{prop}
\begin{proof}
 Suppose that $\tilde{H}_k(G) \neq 0$ for some $k \in \Z$. Then diagonality implies that the components of the chain complex $\cdots \too E^{k, 2}_{k} \too E^{k-2, 2}_{k-1} \too \cdots$ are all $0$ except for $E^{k, 2}_{k} = \tilde{H}_k(G)$, whence we have $E^{k, 3}_{k} = E^{k, 2}_{k}$. Inductively, we obtain $E^{k, \infty}_{k} = \tilde{H}_k(G) \neq 0$, which contradicts Lemma \ref{spgr} and Proposition \ref{rknull} (See also Figure \ref{ill2}). This completes the proof.
\end{proof}
\begin{figure}[htbp]
\centering
\begin{tikzpicture}[scale=1.2]
\filldraw[fill=black, draw=black] (0, 0)  node[left] {$\dMH^{0}_0$} ;
\filldraw[fill=black, draw=black] (1, 0)  node[left] {$\dMH^{1}_1$} ;
\filldraw[fill=black, draw=black] (2, 0)  node[left] {$\dMH^{2}_2$} ;
\filldraw[fill=black, draw=black] (3, 0)  node[left] {$\dMH^{3}_3$} ;
\filldraw[fill=black, draw=black] (4, 0)  node[left] {$\dMH^{4}_4$} ;

\filldraw[fill=black, draw=black] (1, -1)  node[left] {$0$} ;
\filldraw[fill=black, draw=black] (2, -1)  node[left] {$0$} ;
\filldraw[fill=black, draw=black] (3, -1)  node[left] {$0$} ;
\filldraw[fill=black, draw=black] (4, -1)  node[left] {$0$} ;

\filldraw[fill=black, draw=black] (2, -2)  node[left] {$0$} ;
\filldraw[fill=black, draw=black] (3, -2)  node[left] {$0$} ;
\filldraw[fill=black, draw=black] (4, -2)  node[left] {$0$} ;

\draw[line width=1pt] (-1,0.5)--(-1,-2.4);
\draw[line width=1pt] (-1.5,-0.5)--(4.5,-0.5);
\draw[->] (2.7,-0.8)--(0.9,-0.2);
\draw[->] (3.7,-1.8)--(0.9,-0.2);
\end{tikzpicture}
\label{graphty}
\caption{This illustrates vanishing of $\dMH^1_1(G)$.}
\label{ill2}
\end{figure}
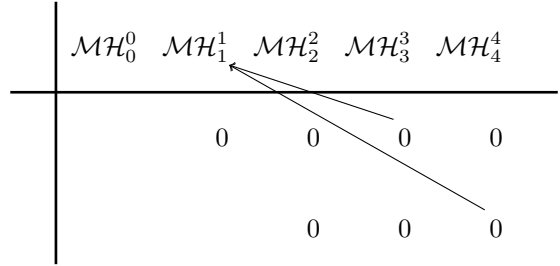
\begin{eg}
Example \ref{diageg} and Proposition \ref{diagvanish} require $\tilde{H}_n(\overrightarrow{T})$ and $\tilde{H}_n(\overrightarrow{K_n})$ to be $0$. Indeed, this is true by Examples \ref{egcontr2}, \ref{egcontr}, and Proposition \ref{htpyequiv}.
\end{eg}
\begin{cor}\label{exact}
If a digraph $G$ is diagonal and has finite diameter, then the sequence
\[
\cdots \too \MH^{\ell}_\ell(G) \too \MH^{\ell-1}_{\ell-1}(G) \too \cdots \too \MH^{0}_{0}(G) \too \Z \too 0
\]
is exact.
\end{cor}
\begin{rem}[M. Yoshinaga]\label{yoshinaga}
 Since the magnitude $\chi_g(q)$ of an undirected diagonal graph $g$ is equal to $\sum_{i=0}^{\infty} (-1)^i\rk \MH^i_i(\overrightarrow{g})q^i$ (\cite{HW} Theorem 2.8), the above Corollary \ref{exact} implies that $\chi_g(q) = \sum_{i=0}^\infty (-1)^i(\rk \del'^i_i + \rk \del'^{i+1}_{i+1})q^i = 1 + (1-q)\sum_{i=0}^{\infty} (-1)^i\rk \del'^{i+1}_{i+1}q^i$. Although it does not necessarily converge as an infinite series, we informally obtain that the limit of magnitude as $q \to 1$ is $1$. Hence Corollary \ref{exact} can be considered as a categorification of the property $\lim_{q \to 1}\chi_g(q) = 1$ that is verified for many graphs (\cite{L1} Example 3.12).
\end{rem}
The following example shows that even if a directed graph $G$ is contractible, it is not necessarily diagonal. Hence the converse of Proposition \ref{diagvanish} is not true.
\begin{eg}\label{counter}
Consider the undirected graph $g$ in Figure 1 and its subgraph $g'$ consisting of vertices $\{0', 1', 2'\}$ and edges between them. We define a digraph map $f : \overrightarrow{g} \too \overrightarrow{g'}$ by $f(i) = i'$ and $f(i') = i'$ for $0\leq i \leq 2$. We denote the inclusion $\overrightarrow{g'} \too \overrightarrow{g}$ by $\iota$. Then we have $f\circ \iota = {\rm Id}_{\overrightarrow{g'}}$ and $\iota \circ f \simeq {\rm Id}_{\overrightarrow{g}}$ by \textit{Remark} \ref{onestep}. Hence we obtain $\overrightarrow{g} \simeq \overrightarrow{g'} \simeq \{0\}$, whence we have $\tilde{H}_k(\overrightarrow{g}) = 0$ for every $k \in \Z$. On the other hand, the calculation in Example 6.5 of \cite{TY} shows that $g$ is not diagonal.
\end{eg}
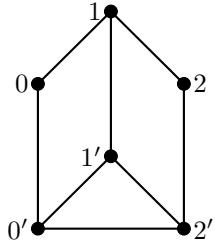
\begin{figure}[htbp]
\centering
\begin{tikzpicture}[scale=1.2]
\filldraw[fill=black, draw=black] (0.8, 2.4) circle (2pt) node[left] {$1$} ;
\filldraw[fill=black, draw=black] (0, 1.6) circle (2pt) node[left] {$0$} ;
\filldraw[fill=black, draw=black] (0, 0) circle (2pt) node[left] {$0'$} ;
\filldraw[fill=black, draw=black] (1.6, 0) circle (2pt) node[right] {$2'$} ;
\filldraw[fill=black, draw=black] (1.6, 1.6) circle (2pt) node[right] {$2$} ;
\filldraw[fill=black, draw=black] (0.8, 0.8) circle (2pt) node[left] {$1'$} ;

\draw[thick] (0.8, 2.4)--(0, 1.6);
\draw[thick] (0, 1.6)--(0, 0);
\draw[thick] (0, 0)--(1.6, 0);
\draw[thick] (1.6, 0)--(1.6, 1.6);
\draw[thick] (1.6, 1.6)--(0.8, 2.4);
\draw[thick] (0.8, 2.4)--(0.8, 0.8);
\draw[thick] (0, 0)--(0.8, 0.8);
\draw[thick] (1.6, 0)--(0.8, 0.8);


\end{tikzpicture}
\label{graphty}
\caption{A contractible and non-diagonal digraph. Each edge is oriented in both directions.}
\end{figure}
\subsection{Non-triviality of $\dMH^{\ell}_k(G)$ for $\ell \neq k$}
\begin{figure}[htbp]
\centering
\begin{tikzpicture}[scale=1.2]
\filldraw[fill=black, draw=black] (0, 0)  node[left] {$E^{0, r}_0$} ;
\filldraw[fill=black, draw=black] (1, 0)  node[left] {$E^{1, r}_1$} ;
\filldraw[fill=black, draw=black] (2, 0)  node[left] {$E^{2, r}_2$} ;
\filldraw[fill=black, draw=black] (3, 0)  node[left] {$E^{3, r}_3$} ;
\filldraw[fill=black, draw=black] (4, 0)  node[left] {$E^{4, r}_4$} ;
\filldraw[fill=black, draw=black] (5, 0)  node[left] {$E^{5, r}_5$} ;
\filldraw[fill=black, draw=black] (6, 0)  node[left] {$E^{6, r}_6$} ;

\filldraw[fill=black, draw=black] (1, -1)  node[left] {$0$} ;
\filldraw[fill=black, draw=black] (2, -1)  node[left] {$E^{2, r}_1$} ;
\filldraw[fill=black, draw=black] (3, -1)  node[left] {$E^{3, r}_2$} ;
\filldraw[fill=black, draw=black] (4, -1)  node[left] {$E^{4, r}_3$} ;
\filldraw[fill=black, draw=black] (5, -1)  node[left] {$E^{5, r}_4$} ;
\filldraw[fill=black, draw=black] (6, -1)  node[left] {$E^{6, r}_5$} ;

\filldraw[fill=black, draw=black] (2, -2)  node[left] {$0$} ;
\filldraw[fill=black, draw=black] (3, -2)  node[left] {$0$} ;
\filldraw[fill=black, draw=black] (4, -2)  node[left] {$E^{4, r}_2$} ;
\filldraw[fill=black, draw=black] (5, -2)  node[left] {$E^{5, r}_3$} ;
\filldraw[fill=black, draw=black] (6, -2)  node[left] {$E^{6, r}_4$} ;

\filldraw[fill=black, draw=black] (3, -3)  node[left] {$0$} ;
\filldraw[fill=black, draw=black] (4, -3)  node[left] {$0$} ;
\filldraw[fill=black, draw=black] (5, -3)  node[left] {$0$} ;
\filldraw[fill=black, draw=black] (6, -3)  node[left] {$E^{6, r}_3$} ;

\draw[line width=1pt] (-1,0.5)--(-1,-3.2);
\draw[line width=1pt] (-1.5,-0.5)--(6,-0.5);
\draw[->] (2.4,-0.8)--(0.9,-0.2);
\draw[->] (3.4,-1.8)--(0.9,-0.2);
\end{tikzpicture}
\label{graphty}
\caption{
The differentials into $E^{1, r}_1$ are all coming from $E^{k, r}_2$ for $k \geq 3$. On the other hand, we have $E^{k, r}_2 = 0$ for $k \geq 5$ by Proposition \ref{diambd} with ${\rm diam}\ \protect\overrightarrow{C_5} = 2$.
}
\label{ill4}
\end{figure}
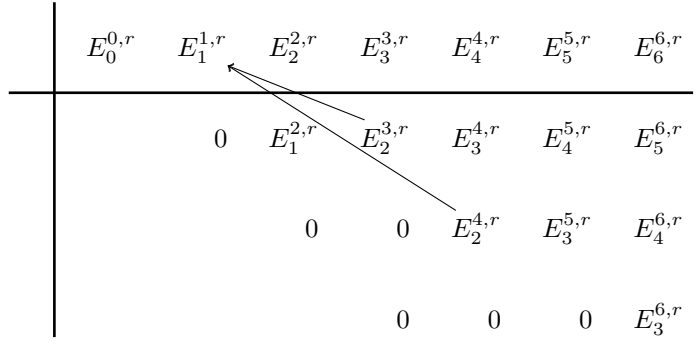

In this subsection, we show that there is a digraph $G$ with $\dMH^{\ell}_k(G) \neq 0$ for some $\ell \neq k$.
\begin{lem}[cf. \cite{GLMY} Example 2.8]\label{5cycle}
 Let $C_5$ be the undirected cycle graph of order $5$. Then we have $\tilde{H}_1(\overrightarrow{C_5}) \neq 0$.
\end{lem}
\begin{proof}
 Let $V(C_5) = \{0, 1, 2, 3, 4\}$ with $\{i, i+1\} \in E(C_5)$ for $0\leq i\leq 3$ and $\{0, 4\} \in E(C_5)$. Then we can check that $(0,1) + (1, 2) + (2, 3) + (3, 4) + (4, 0) \in \Omega^{\rm reg}_1(\overrightarrow{C_5})$ represents a non-zero homology cycle in $\tilde{H}_1(\overrightarrow{C_5})$ as follows. It is immediately checked that 
 \[
 \del^{\rm reg}_1((0,1) + (1, 2) + (2, 3) + (3, 4) + (4, 0)) = 0.
 \]
 Suppose that there is a chain $\alpha \in \Omega^{\rm reg}_2(\overrightarrow{C_5})$ such that 
 \[
 \del^{\rm reg}_2 \alpha = (0,1) + (1, 2) + (2, 3) + (3, 4) + (4, 0).
 \]
 Let $V^+$ and $V^-$ be multi-sets of tuples such that $\alpha = \sum_{x \in V^+}x - \sum_{y \in V^-}y$. We can assume that $V^+\cup V^-$ is minimal among such multi-sets and such $\alpha$'s. This assumption implies that there are no subsets $U^+ \subset V^+$ and $U^- \subset V^-$ such that $\del^{\rm reg}_2(\sum_{x \in U^+}x - \sum_{y \in U^-}y) = 0$. Since the tuple $(0, 1)$ appears in the RHS of the above equation, the tuple $(0, 1, 0)$ or $(1, 0, 1)$ must be in $V^+$ (other possible allowed $2$-paths $(4, 0, 1)$ and $(0, 1, 2)$ are not contained in $\Omega^{\rm reg}_2(\overrightarrow{C_5})$). In both case, the tuple $(0, 1, 0)$ or $(1, 0, 1)$ must be in $V^-$ to cancel the tuple $(1, 0)$. Hence we obtain one-element subsets $U^+ \subset V^+$ and $U^- \subset V^-$ such that $\del^{\rm reg}_2(\sum_{x \in U^+}x - \sum_{y \in U^-}y) = 0$, which contradicts to the minimality assumption. This completes the proof.
\end{proof}
\begin{prop}\label{nontrivial}
 We have $\dMH^{\ell}_k(\overrightarrow{C_5}) \neq 0$ for some $\ell \neq k$. More specifically, we have $\dMH^3_2(\overrightarrow{C_5}) \neq 0$ or $\dMH^4_2(\overrightarrow{C_5}) \neq 0$.
\end{prop}
\begin{proof}
Suppose that $\dMH^3_2(\overrightarrow{C_5}) = \dMH^4_2(\overrightarrow{C_5}) = 0$. Then we have $E^{1, \infty}_1(\overrightarrow{C_5}) = E^{1, 2}_1(\overrightarrow{C_5}) = \tilde{H}_1(\overrightarrow{C_5})$ by Proposition \ref{diambd} and by the similar argument as the proof of Proposition \ref{diagvanish} (See also Figure \ref{ill4}). Since we have $\tilde{H}_1(\overrightarrow{C_5}) \neq 0$ by Lemma \ref{5cycle}, it contradicts Lemma \ref{spgr} and Proposition \ref{rknull}. This completes the proof.
\end{proof}
\subsection{Applications for path homology}
In this subsection, we show some properties of reduced path homology by using following fact proved in \cite{AHK}. Note that the tuple $(x_0, x_1, x_0, \dots)$ obtained by $(k+1)$ times repeating distinct vertices $x_0$ and $x_1$ with $\{x_0, x_1\} \in E(g)$ is a homology cycle of $\MC^{k}_k(\overrightarrow{g})$ for any undirected graph $g$.
\begin{lem}[\cite{AHK} Corollary 1.2]\label{ahk}
If the girth of an undirected graph $g$ is $\geq 5$, then we have 
\[
\MH^{k}_k(\overrightarrow{g}) \cong \Z\langle(x_0, x_1, x_0, \dots) \in V(\overrightarrow{g})^{k+1} \mid (x_0, x_1) \in E(\overrightarrow{g})\rangle,
\]
for $k \geq 1$, and $\MH^{0}_0(\overrightarrow{g}) \cong \Z\langle V(\overrightarrow{g}) \rangle$.
\end{lem}
\begin{prop}\label{girth5}
  If the girth of an undirected graph $g$ is $ \geq 5$, then we have $\tilde{H}_k(\overrightarrow{g}) = 0$ for $k\geq 2$.
\end{prop}
\begin{proof}
 By Proposition \ref{magspec}, the path homologies $\tilde{H}_k(\overrightarrow{g})$ are the homologies of the chain complex 
 \[
 \cdots \too \MH^{k}_k(\overrightarrow{g}) \too \MH^{k-1}_{k-1}(\overrightarrow{g}) \too \cdots \too \MH^{0}_0(\overrightarrow{g}) \too \Z \too 0.
 \]
 By Lemma \ref{ahk}, the differentials of the above chain complex are determined by 
 \[
 \del'^{2m+1}_{2m+1}(x_0, x_1,  \dots, x_1) = (x_1, x_0, \dots, x_1) - (x_0, x_1, \dots, x_0)
 \]
 and 
 \[
 \del'^{2m}_{2m}(x_0, x_1,  \dots, x_0) = (x_1, x_0, \dots, x_0) + (x_0, x_1, \dots, x_1),
 \]
 for $m \geq 0$. Then we can check that 
 \[
 {\mathrm Ker}\ \del'^{2m+1}_{2m+1} = \Z\langle (x_0, x_1, \dots, x_1) + (x_1, x_0, \dots, x_0) \rangle = \Z\langle \del'^{2m}_{2m}(x_0, x_1,  \dots, x_0) \rangle
 \]
 and 
 \[
 {\mathrm Ker}\ \del'^{2m}_{2m} = \Z\langle (x_0, x_1, \dots, x_0) - (x_1, x_0, \dots, x_1) \rangle = \Z\langle \del'^{2m+1}_{2m+1}(x_1, x_0,  \dots, x_0) \rangle,
 \]
 for $m \geq 1$. From these equalities, we obtain $\tilde{H}_k(\overrightarrow{g}) = 0$ for $k\geq 2$. This completes the proof.
 
\end{proof}
\begin{prop}\label{girth52}
  If the girth of an undirected graph $g$ is $ \geq 5$ and $< \infty$, then  we have $\tilde{H}_1(\overrightarrow{g}) \neq 0$.
\end{prop}

\begin{proof}
 By the assumption, there exists an embedding $\iota : C_n \too g$, where $C_n$ is an undirected cycle graph with $n \geq 5$. We denote the embedded vertices by $x_0, \dots, x_{n-1}$ with $\{x_i, x_{i+1}\}, \{x_0, x_{n-1}\} \in E(G)$ and $x_i \neq x_j$ for $0 \leq i \leq n-2$ and $j \neq i$. Then we have 
 \[
 (x_0, x_1) + \cdots + (x_{n-2}, x_{n-1}) + (x_{n-1}, x_0) \in {\mathrm Ker}\ \del'^{1}_1.
 \]
 Since we have $\MH^{2}_2(\overrightarrow{g}) \cong \Z\langle(x_0, x_1, x_0) \in V(\overrightarrow{g})^{3} \mid (x_0, x_1) \in E(\overrightarrow{g})\rangle$, we can check that 
 \[
 (x_0, x_1) + \cdots + (x_{n-2}, x_{n-1}) + (x_{n-1}, x_0) \not\in {\mathrm Im}\ \del'^{2}_2
 \]
 similarly to the proof of Lemma \ref{5cycle}. Hence we obtain $0 \neq [(x_0, x_1) + \cdots + (x_{n-2}, x_{n-1}) + (x_{n-1}, x_0)] \in \tilde{H}_1(\overrightarrow{g})$. This completes the proof.
\end{proof}
\begin{rem}\label{errortum}
 In Example 6.10 (i) of \cite{GLMY}, they state that \[
 \tilde{H}_k(\overrightarrow{C_p}) = \begin{cases} \Z & k = p, \\ 0 & \text{otherwise},\end{cases}
 \]
 for $p \geq 5$ without a proof. From our results, this should be corrected to \[
 \tilde{H}_k(\overrightarrow{C_p}) = \begin{cases} \Z & k = 1, \\ 0 & \text{otherwise}\end{cases}.\]
\end{rem}
\begin{rem}
 In \cite{AHK}, they show the following more detailed result than Lemma \ref{ahk}.
 \end{rem}
 \begin{prop}[\cite{AHK} Corollary 1.4]
 Let $\ell \geq 1$ and $i \geq 0$. If the girth of an undirected graph $g$ is $\geq 2i+5$, then 
 \[
 \MH^{\ell}_{\ell - j} (\overrightarrow{g}) = \begin{cases} \Z^{\# E(g)} & j = 0 \\ 0 & 1 \leq j \leq i \end{cases}.\]
 \end{prop}
 This proposition implies the following by the same argument in the proof of Propositions \ref{girth5} and \ref{girth52}, and by Proposition \ref{magspec}.
\begin{prop}
 Let $\ell \geq 1$ and $i \geq 0$. If the girth of an undirected graph $g$ is $\geq 2i+5$, then $\dMH^{\ell}_{\ell - j} (\overrightarrow{g}) = 0$ for $j = 0$ with $\ell \geq 2$ or $1 \leq j \leq i$, and $\dMH^{1}_{1} (\overrightarrow{g}) \neq 0$.
 \end{prop}



\begin{thebibliography}{99}
\bibitem{A} Y.Asao, \textit{Magnitude homology of geodesic metric spaces with an upper curvature bound}, Alg. Geom. Topol. 21 (2021) 647--664.
\bibitem{AHK} Y. Asao, Y. Hiraoka and S. Kanazawa, \textit{Girth, magnitude homology, and phase transition of diagonality}, preprint, arXiv:2101.09044, 2021.
\bibitem{Gom}K. Gomi, \textit{Magnitude homology of geodesic space}, preprint, 
arXiv:1902.07044, 2019.
\bibitem{Gu}Y. Gu, \textit{Graph magnitude homology via algebraic Morse theory}, preprint, 
arXiv:1809.07240, 2018.
\bibitem{GJMY}A. Grigor'yan, R. Jimenez, Y. Muranov and S.-T. Yau, \textit{On the path homology theory of digraphs and Eilenberg-Steenrod axioms}, Homology Homotopy Appl. \textbf{20} (2018), no. 2, 179--205. \bibitem{GJMY2}A. Grigor'yan, R. Jimenez, Y. Muranov and S.-T. Yau, \textit{Homology of path complexes and hypergraphs}, Topology and its Appl. \textbf{267} (2019), no. 2, 106877.
\bibitem{GLMY}A. Grigor'yan, Y. Lin, Y. Muranov and  S.-T. Yau, \textit{Homotopy theory for digraphs}, Pure Appl. Math. Q. \textbf{10} (2014), no. 4, 619--674. 
\bibitem{GMY}A. Grigor'yan, Y. Muranov and  S.-T. Yau, \textit{On a cohomology of digraphs and Hochschild cohomology}, J. Homotopy and related structures \textbf{11} (2016), 209--230.
\bibitem{GMVY}A. Grigor'yan, Y. Muranov, V. Vershinin and  S.-T. Yau, \textit{Path homology theory of multigraphs and quivers }, Forum Mathematicum \textbf{30} (2018), no. 5, 1319--1337.
\bibitem{HW}
R. Hepworth and S. Willerton,
\textit{Categorifying the magnitude of a graph}, Homology Homotopy Appl. \textbf{19} (2017), 31--60.
\bibitem{YK}R. Kaneta and  M. Yoshinaga,
\textit{Magnitude homology of metric spaces and order complexes}, Bulletin of the London Mathematical Society \textbf{53}(3) (2021), 893--905.
\bibitem{L1} T. Leinster, \textit{The magnitude of a graph}
, Mathematical Proceedings of the Cambridge Philosophical Society \textbf{166} (2019), 247--264.
\bibitem{LS}T. Leinster and  M. Shulman, \textit{Magnitude homology of enriched categories and metric spaces}, Alg. Geom. Topol. \textbf{21} (2021), 2175--2221.
arXiv:1912.13483, 2019.
\bibitem{TY}Y. Tajima and M. Yoshinaga, \textit{Magnitude homology of graphs and discrete Morse theory on Asao–Izumihara complexes}, preprint, 
arXiv:2110.02458, 2021.
\bibitem{Wei}C.A. Weibel, \textit{An introduction to homological algebra}, Cambridge Studies in Advanced Mathematics, 38. Cambridge University Press, Cambridge, 1994.
\end{thebibliography}
\end{document}